\documentclass[12pt]{article}
\usepackage{amsmath,amsthm,amsfonts}

\usepackage{subfigure}
\usepackage[mathscr]{eucal}
\usepackage{graphics,epsfig}
\usepackage{longtable}
\usepackage{algorithmic}
\usepackage[ruled]{algorithm}
\usepackage{multirow}
\usepackage{rotating}
\usepackage{verbatim}

\usepackage{color}

\numberwithin{equation}{section} \numberwithin{figure}{section}
\numberwithin{table}{section}

\def\Box{\sharp}

\def\bb{\begin{equation}} \def\ee{\end{equation}}

\def\beqn{\begin{eqnarray}}  \def\eqn{\end{eqnarray}}

\def\beqnx{\begin{eqnarray*}} \def\eqnx{\end{eqnarray*}}

\def\bn{\begin{enumerate}} \def\en{\end{enumerate}}

\def\bd{\begin{description}} \def\ed{\end{description}}

\begin{document}

\begin{center}
{\bf Adaptive Uzawa algorithm for nonsymmetric generalized saddle
point problem}
\end{center}

\bigskip
\centerline{Hailun ~Shen\footnote{School of Mathematics and
Statistics, Wuhan University, Wuhan 430072, P. R. China. ({\tt
2012202010045@whu.edu.cn}).}
\quad \quad Hua ~Xiang\footnote{School of Mathematics and
Statistics, Wuhan University, Wuhan, China. ({\tt
hxiang@whu.edu.cn}). The work was completed while Hua Xiang was
visiting The Chinese University of Hong Kong in 2013 under the
support of Hong Kong RGC grant (Project 405110).}  }

\bigskip
\begin{abstract}
In this paper, we extend the inexact Uzawa algorithm in [Q. Hu, J.
Zou, SIAM J. Matrix Anal., 23(2001), pp. 317-338] to the
nonsymmetric generalized saddle point problem. The techniques used
here are similar to those in [Bramble \emph{et al}, Math. Comput.
69(1999), pp. 667-689], where the convergence of Uzawa type
algorithm for solving nonsymmetric case depends on the spectrum of
the preconditioners involved. The main contributions of this paper
focus on a new linear Uzawa type algorithm for nonsymmetric
generalized saddle point problems  and its convergence. This new
algorithm can always converge without any prior estimate on the
spectrum of two preconditioned subsystems involved,  which may not
be easy to achieve in applications. Numerical results of the
algorithm on the Navier-Stokes problem are also presented.\\
\end{abstract}

%%%%%%%%%%%%%%%%%%%%%%%%%%%%%%%%%%%%%%%
%
\section{Introduction}
%
%%%%%%%%%%%%%%%%%%%%%%%%%%%%%%%%%%%%%%%

Let $H_1$ and $H_2$ be finite dimensional Hilbert spaces with inner
products denoted by $\langle \cdot  , \cdot \rangle $
(cf.\cite{BramblePasciakVassilev_MC99}). We consider to solve the
following system
\begin{equation}\label{eq1.0}
 \begin{bmatrix} A & B \\
 B^{T}& -D \end{bmatrix}\begin{bmatrix} x\\ y\end{bmatrix}=\begin{bmatrix} f\\g\end{bmatrix},
\end{equation}
where $A$ is an $n \times n $ nonsymmetric matrix, $B$ is an ${n
\times m}$ matrix with ${m \leq n}$, and $D$ is a symmetric
semi-positive matrix. We shall assume that the Schur complement
matrix
 \begin{equation*}
 S=B^TA^{-1}B+D
 \end{equation*}
 is nonsingular.

The system \eqref{eq1.0} arises from many areas of computational
sciences and engineerings, for example, in certain finite element
and finite difference discretization of Navier-Stokes equations,
Oseen equations, and mixed finite element discretization of second
order convection-diffusion problems (cf. \cite{M. B MAA05
,VGPA_MAA81,  O. A. K MAA82,  P.K01, R. T77, P. S. V, H. L10}). For
the saddle point problem, there exist many algorithms, for example,
the Krylov iteration methods with block diagonal, triangular block
or constraint preconditioners (see \cite{M. B MAA05 } and the
references therein). The Uzawa type algorithms applied to
nonsymmetric saddle point problems are of great interest because
they are simple, efficient, and have minimal computer memory
requirements. They can be applied to the solution of difficult
practical problems such as the Navier-Stokes equation. Many
algorithms are applied to the system \eqref{eq1.0}  when $A$ is a
symmetric positive definite matrix (see \cite{B.W08,R. B.
H_MAA90,BramblePasciakVassilev_MC97,ZHCAO_MAA03,HuZou_MAA01,
HuZou_MAA02, HuZou_MAA06, J.L.Z.Z_MAA10,Z. T. A. S_MAA98} and the
references therein). Bramble, Pasciak and Vassilev
\cite{BramblePasciakVassilev_MC99} investigated the convergence of
Uzawa method on nonsymmetric  saddle point problem.
Cao \cite{ZHCao_MAA04} considered generalized saddle point problems
with $D\neq 0$ and the acceleration of the convergence of the
inexact Uzawa algorithms, together with a new nonlinear Uzawa type
algorithm.
But nearly all existing preconditioned Uzawa algorithms for
nonsymmetric case do not adopt self-updating relaxation parameters,
and converge only under some proper scalings of the preconditioners
$A_s$ and $B^TA_s^{-1}B+D $, where $A_s$ is the symmetric part of
$A$.
Hu and Zou \cite{HuZou_MAA01} suggested Uzawa type algorithm for
symmetric saddle point problems with variable relaxation parameters.
But few studies on the convergence analysis of preconditioned Uzawa
method can be found for  nonsymmetric saddle point problems with
relaxation parameters. In this paper, we combine the techniques in
\cite{BramblePasciakVassilev_MC99} and \cite{HuZou_MAA01}. We extend
the Uzawa algorithm with variable parameter for symmetric saddle
point problem in  \cite{HuZou_MAA01} to the nonsymmetric case, and
also modified the Uzawa algorithm for nonsymmetric saddle point
problem in \cite{BramblePasciakVassilev_MC99} with variable
parameter.

\indent Throughout this paper we assume that $A$ has a positive
definite symmetric part. The symmetric part $A_s$ of the operator
$A$ is defined by
\begin{equation}
A_s=\frac{1}{2}(A+A^{T}). \\ \end{equation}
\indent We assume that $A_s$ is positive definite and satisfies\\
\begin{equation}\label{eq1.1}
\langle Ax,y\rangle\leq\alpha\langle A_sx,x\rangle^{1/2}\langle
A_sy,y\rangle^{1/2},
\end{equation}
for all $x,y\in H_1$.
Under this assumption, the system \eqref{eq1.0} is solvable if and only if LBB condition is assumed to hold for the pair of spaces $H_1$ and $H_2$, i.e.,\\
\begin{equation}
\sup_{u\in{H_1}}\frac{\langle v,Bu\rangle^2}{\langle
A_su,u\rangle}\geq c_0\Vert v\Vert^2, ~ \forall   \ v\in{H_2},
\end{equation}
for some positive number $c_0$. Here $\Vert.\Vert$ denote the norm
in the space $H_2$ (or $H_1$) corresponding to the inner product
$\langle\cdot,\cdot\rangle$. See Theorem 2.1 in
 \cite{BramblePasciakVassilev_MC99}.

 Our algorithms are motivated by  Uzawa iteration with variable relaxation parameters for symmetric saddle point problems with $D=0$ in \cite{HuZou_MAA01}, which can be defined as follows.\\
  {\bf Algorithm 1.1} (Hu-Zou \cite{HuZou_MAA01}). Given $x_0\in
{H_1}$ and $y_0\in {H_2}$, the sequence$\{(x_i,y_i)\}$ is defined,
for ${i=0,1,2,...,}$ by
\begin{eqnarray*}
x_{i+1} & = & x_i+\omega_i\widehat A^{-1}(f-(Ax_i+By_i)),\\
y_{i+1}&=&y_i+\theta\tau_i\widehat{C}^{-1}(B^{T}x_{i+1}-g),
\end{eqnarray*}
where $A$ is symmetric positive definite. The relaxation parameter
$\tau_i $ is determined such that the norm
\begin{equation*}
\Vert \tau_i\widehat C^{-1}g_i-C^{-1}g_i\Vert_{C}^2
 \end{equation*}
 is minimized, where $g_i=B^Tx_{i+1}-g, C=B^TA^{-1}B$, and $\widehat C $ is the preconditioner for $C$. Then we choose
 \begin{equation*}
 \tau_i= \left\{
   \begin{array}{cc}
  \frac{\langle g_i , \widehat{C}^{-1}g_i\rangle}{\langle B^TA^{-1}B\widehat{C}^{-1}g_i , \widehat{C}^{-1}g_i\rangle},&g_i\neq0,  \\
   1,&g_i=0. \\
   \end{array}
   \right.
  \end{equation*}
  The relaxation parameter $\omega_i$ is determined such that the norm
  \begin{equation*}
\Vert A^{-1}f_i-\omega_i\widehat A^{-1}f_i\Vert_{A}^2
 \end{equation*}
  is minimized, where $f_i=f-Ax_i-By_i$, $r_i=\widehat A^{-1}f_i,$ $\widehat A$ is the preconditioner for
  $A$, and we set
\begin{equation*}
 \omega_i= \left\{
   \begin{array}{cc}
  \frac{\langle f_i,r_i\rangle}{\langle Ar_i,r_i\rangle},&f_i\neq0,  \\
   1,&f_i=0. \\
   \end{array}
   \right.
  \end{equation*}
\\ \indent We are concerned about  whether this algorithm can be applied to nonsymmetric case, which will be discussed later.
For the nonsymmetric matrix $A$ and $D=0$, Bramble, Pasciak and
Vassilev presented the linear inexact Uzawa algorithm in
\cite{BramblePasciakVassilev_MC99} as
follows.\\
{\bf Algorithm 1.2} (Bramble-Pasciak-Vassilev
\cite{BramblePasciakVassilev_MC99}). Given $x_0\in {H_1}$ and
$y_0\in {H_2}$, the sequence$\{(x_i,y_i)\}$ is defined, for
${i=0,1,2,...,}$ by
\begin{eqnarray*} x_{i+1} & = & x_i+\delta A_0^{-1}(f-(Ax_i+By_i)),
\\
y_{i+1}&=&y_i+\tau Q_B^{-1}(B^Tx_{i+1}-g).
\end{eqnarray*}
Here $\tau$ and $\delta$ are positive constant parameters, $A_0$ and
$Q_B$ are the preconditioners for $A_s$ and $B^TA_s^{-1}B$,
respectively, and satisfying
\begin{equation}\label{eq1.11}
\langle A_0v,v\rangle\leq \langle A_sv,v\rangle \leq \kappa_0\langle
A_0v,v\rangle,
\end{equation}
for all $v\in{H_1}$, and
\begin{equation}\label{eq1.12}
\gamma\langle Q_Bw,w\rangle\leq\langle
B^TA_s^{-1}Bw,w\rangle\leq\langle Q_Bw,w\rangle,
\end{equation}
for all $w\in {H_2}$, where $\gamma \in [0,1]$.

 \indent The
inequalities \eqref{eq1.11} and \eqref{eq1.12} respectively imply
scaling of $A_0$ and $Q_B$. That is to say, this algorithm is
convergent only under the proper scaling of the preconditioners,
which is not be easy to achieve in applications. Therefore, we
suggest an algorithm in this paper to overcome the limitations
above. Moreover, our algorithm is applied to  generalized saddle
point problem \eqref{eq1.0} for $D\neq 0$.

The paper is organized as follows. In section 2 we analyze an exact
Uzawa algorithm for solving \eqref{eq1.0}. In section 3 we define
and analyze a linear one-step Uzawa type algorithm. Section 4
provides the results of numerical experiments.

%\newpage \vfil  \vfil
For the sake of clarity, we list the main notations used later.\\
% \makebox[10cm]{\hrulefill} \\
% $S=B^TA^{-1}B+D$, the exact Schur complement of \eqref{eq1.0} \\
% $A_s,S_s$ are spd parts of $A$ and $S$, respectively\\
% $H=B^TA_s^{-1}B+D$\\
% $\widehat S ,A_0 $ are the spd preconditioners of the matrices $H$ and $A_s$, respectively \\
% $\kappa_1=cond(\widehat S^{-1}H), \beta1=\frac{\kappa_1-1}{\kappa_1+1}$\\
% $\kappa_2=cond(\widehat S^{-1}S_s), \beta2=\frac{\kappa_2-1}{\kappa_2+1}$\\
% $\kappa_3=cond(\widehat S^{-1}M), \beta3=\frac{\kappa_3-1}{\kappa_3+1}$, where  $M=B^TA_0^{-1}B+D$\\
% \makebox[10cm]{\hrulefill}
%
 \begin{tabular}{l} \hline
 $S=B^TA^{-1}B+D$, the exact Schur complement of \eqref{eq1.0} \\
 $A_s,S_s$ are spd parts of $A$ and $S$, respectively\\
 $H=B^TA_s^{-1}B+D$\\
 $\widehat S ,A_0 $ are the spd preconditioners of the matrices $H$ and $A_s$, respectively \\
 $\kappa_1=cond(\widehat S^{-1}H), \beta1=\frac{\kappa_1-1}{\kappa_1+1}$\\
 $\kappa_2=cond(\widehat S^{-1}S_s), \beta2=\frac{\kappa_2-1}{\kappa_2+1}$\\
 $\kappa_3=cond(\widehat S^{-1}M), \beta3=\frac{\kappa_3-1}{\kappa_3+1}$, where $M=B^TA_0^{-1}B+D$\\ \hline
 \end{tabular}

\section{Analysis of the preconditioned exact Uzawa algorithm with relaxation parameter}

In this section, we first give an exact Uzawa algorithm for \eqref{eq1.0} with the nonsymmetric matrix $A$, then analyze the convergence of this algorithm. The preconditioned variant of the exact Uzawa algorithm with relaxation parameter is defined as follows.\\
 {\bf Algorithm 2.1} Given $x_0\in
{H_1}$ and $y_0\in {H_2}$, the sequence$\{(x_i,y_i)\}$ is defined,
for ${i=0,1,2,...}$, by
\begin{eqnarray*}
x_{i+1} & = & x_i+A^{-1}(f-(Ax_i+By_i)),\\
y_{i+1}&=&y_i+\theta\tau_i\widehat{S}^{-1}(B^{T}x_{i+1}-Dy_i-g).
\end{eqnarray*}
And the relaxation parameter $\tau_i $ is determined such that the
norm
\begin{equation}\label{eq2.01}
\Vert \tau_i\widehat S^{-1}g_i-H^{-1}g_i\Vert_{H}^2.
 \end{equation}
 is minimized, where $g_i=B^Tx_{i+1}-Dy_i-g, H=B^TA_s^{-1}B+D,$ then\\
 \begin{equation}\label{eq0.00}
 \tau_i= \left\{
   \begin{array}{cc}
  \frac{\langle g_i,\widehat S^{-1}g_i\rangle}{\langle(B^TA_s^{-1}B+D)\widehat S^{-1}g_i,\widehat S^{-1}g_i\rangle},&g_i\neq0,  \\
   1,&g_i=0. \\
   \end{array}
   \right.
  \end{equation}
%\begin{equation}
%\tau_i=\left\{
%\begin{array}{c}\frac{(g_i,\hat{C}^{-1})}{(A^{-1}B\hat{C}^{-1},B\hat{C}^{-1})},g_i\leq0.\\1,g_i=0
%\end {array} \right
%\end{equation}

The relaxation parameter  above can be computed effectively, similar
to the evaluation of the  iteration parameter in the conjugate
gradient method. In this paper, we follow \cite{HuZou_MAA01} to
evaluate the parameter $\tau_i$, but here $A$ in Algorithm 2.1 is
nonsymmetric. So we make some small modification in the choose of
$\tau_i$. It will be shown that our algorithm always converges for
general preconditioner $\widehat{S}$, while the convergence of most
existing Uzawa-type algorithms for solving nonsymmetric saddle point
is guaranteed only under certain conditions on the extreme
eigenvalues of the preconditioned matrix $\widehat{S}^{-1}H$.

Define the iteration errors of the above method by %$e_i^x = x-x_i$ and $e_i^y = y-y_i$.
\begin{eqnarray*}
e_i^x &=& x-x_i, \\
e_i^y &=& y-y_i.
\end{eqnarray*}
%be generated by the above method.
We can derive
\begin{equation*}
e_{i+1}^y=(I-\theta\tau_i\widehat{S}^{-1}(B^{T}A^{-1}B+D))e_i^y=(I-\theta\tau_i\widehat{S}^{-1}S)e_i^y.
\end{equation*}
Therefore, the convergence of Algorithm 2.1 is governed by the
properties of the operator $(I-\theta\tau_i\widehat{S}^{-1}S)$.

In order to prove the convergence of Algorithm 2.1, we need the following lemma.\\
{\bf Lemma 2.1}  Suppose that $A$ is an invertible linear operator
with positive definite symmetric part $A_s$ and satisfies
\eqref{eq1.1}. Then $(A^{-1})_s$ is positive definite and satisfies
\begin{equation}
\langle(A^{-1})_sw,w\rangle\leq\langle
A_s^{-1}w,w\rangle\leq\alpha^2\langle(A^{-1})_sw,w\rangle, ~ \forall
\ w \in H_1.
\end{equation}
See Lemma 2.1 in \cite{BramblePasciakVassilev_MC99}.

{\bf Lemma 2.2} For any natural number $i$, there is a symmetric and positive definite ${m \times m}$ matrix ${ G_i }$ such that\\
 (i) $G_i^{-1}g_i=\theta\tau_i\widehat{S}^{-1}g_i $ with $g_i=B^Tx_{i+1}-Dy_i-g $ as defined in Algorithm 2.1;\\
 (ii) All eigenvalues of the matrix $ G_i^{-1}H $ lie in the interval $[\theta(1-\beta_1),\theta(1+\beta_1)]$.\\
  Proof: By the definition of the parameter $\tau_i$ we have
 \begin{equation*}
 \begin {split}
 &||\tau_i\widehat{S}^{-1}g_i-H^{-1}g_i||_H^2\\ =& ||H^{-1}g_i||_H^2-2\tau_i\langle g_i,\widehat{S}^{-1}g_i\rangle +\tau_i\frac{\langle g_i,\widehat S^{-1}g_i\rangle}{\langle(B^TA_s^{-1}B+D)\widehat{S}^{-1}g_i,\widehat{S}^{-1}g_i\rangle}\Vert\widehat{S}^{-1}g_i\Vert_H^2\\
 = & ||H^{-1}g_i||_H^2-2\tau_i\langle g_i,\widehat{S}^{-1}g_i\rangle +\tau_i\langle g_i,\widehat{S}^{-1}g_i\rangle\\=&(1-\tau_i\frac{ \langle g_i,\widehat{S}^{-1}g_i\rangle}{\langle
 g_i,H^{-1}g_i\rangle})||H^{-1}g_i||_H^2.
\end {split}
\end{equation*}
It follows from the well-known Kantororich inequality that
\begin{equation}
\frac{\langle v,v\rangle\langle v,v\rangle}{\langle
Gv,v\rangle\langle
G^{-1}v,v\rangle}\geq\frac{4\lambda_1\lambda_2}{(\lambda_1+\lambda_2)^2},
~ \forall \ v\in \mathbb{R}^l,
\end{equation}
where $\lambda_1$ and $\lambda_2$ are the smallest and largest eigenvalues of the $l\times l$ symmetric positive matrix $G$.\\
Then from the definition of $\tau_i$ in \eqref{eq0.00} and the
well-known inequality above, we obtain
\begin{equation*}
\begin {split}
&\tau_i\frac{ \langle g_i,\widehat{S}^{-1}g_i\rangle}{\langle g_i,
 H^{-1}g_i\rangle}\\ =& \frac{\langle\widehat{S}^{-1/2}g_i,\widehat{S}^{-1/2}g_i\rangle^2}{\langle\widehat{S}^{-1/2}H\widehat{S}^{-1/2}(\widehat{S}^{-1/2}g_i),\widehat{S}^{-1/2}g_i\rangle\langle\widehat{S}^{1/2}H^{-1}\widehat{S}^{1/2}(\widehat{S}^{-1/2}g_i),\widehat{S}^{-1/2}g_i\rangle}\\
 \geq & \frac{4\lambda_1^{'}\lambda_2^{'}}{(\lambda_1^{'}+\lambda_2^{'})^2}
 =\frac{4\kappa_1}{(1+\kappa_1)^2},
\end {split}\end{equation*}\\
where $\lambda_1^{'}$ and $\lambda_2^{'}$ are the minimal and maximal eigenvalues of the matrix $\widehat{S}^{-1/2}H\widehat{S}^{-1/2}$,\\respectively. Hence we obtain\\
\begin{equation*}||\tau_i\widehat{S}^{-1}g_i-H^{-1}g_i||_H\leq(1-\frac{4\kappa_1}{(1+\kappa_1)^2})||H^{-1}g_i||_H=\beta_1||H^{-1}g_i||_H.\end{equation*}
 It is clear that\ $\beta_1<1$.
 This implies the existence of a symmetric positive definite ${m \times m}$ matrix $\widehat G_{i}$ such that
\begin{equation*}\widehat G_{i}^{-1}g_i=\tau_i\widehat{S}^{-1}g_i ,\end{equation*}
and
 \begin{equation*}||I-H^{1/2}\widehat G_{i}^{-1}H^{1/2}||\leq\beta_1.\end{equation*}
 See Lemma 9 in \cite{R. B. H_MAA90}, then the existence of such a matrix $\widehat G_{i}$ is proved.

 Now set $G_{i}^{-1}=\theta \widehat G_{i}^{-1}$,  and then we have
 \begin{equation*} G_{i}^{-1}g_i=\theta\tau_i\widehat{S}^{-1}g_i. \end{equation*}
 And we also know that all eigenvector of the matrix $H^{1/2}G_{i}^{-1}H^{1/2}$ lie in the interval $[\theta(1-\beta_1),\theta(1+\beta_1)]$,
 which yields the desired eigenvalue bounds.~~~~~~~~~~~~~~~~~~~~~~~~~~~~~~~~~~~~~~~~~~~~~~~~~~~~~~~~~~~~~~~~~~~~~~~~~~$\Box$

{\bf Theorem 2.1}  Suppose that $A$ is invertible with positive definite symmetric part $A_s$ which satisfies \eqref{eq1.1}. Suppose also that $A_s$ satisfies LBB condition. Then we have\\
 \begin{equation*} ||(I-\theta\tau_i\widehat{S}^{-1}(B^TA^{-1}B+D))u||^{2}_{G_{i}}\leq(1-\frac{\theta(1-\beta_1)}{\alpha^2})||u||^2_{G_{i}}, ~\forall \ u\in{H_2},
\end{equation*}
when $0<\theta<\frac{1-\beta_1}{\alpha^2(1+\beta_1)^2}$.\\
Proof:   We set $\Gamma=B^{T}A^{-1}B$, then according to Lemma 2.2,
 \begin{eqnarray}\label{eq2.1}
\begin{split}
&||(I-\theta\tau_i\widehat{S}^{-1}(B^{T}A^{-1}B+D))u||^{2}_{G_{i}}\\=&||(I-G_{i}^{-1}(B^{T}A^{-1}B+D))u||^{2}_{G_{i}}\\=&||u||^{2}_{G_{i}}-2\langle(\Gamma
+D) u,u\rangle+\langle(\Gamma+D) u,G_{i}^{-1}(\Gamma +D)u\rangle.
\end{split}
\end{eqnarray}
In addition, according to Cauchy-Schwarz inequality and Lemma 2.1,
we have
\begin{eqnarray}\label{eq2.02}
\begin{split}
\langle A^{-1}v,w\rangle&=\langle
A_s^{1/2}A^{-1}v,A_s^{-1/2}w\rangle\\&\leq\langle
A_s^{1/2}A^{-1}v,A_s^{1/2}A^{-1}v\rangle^{1/2} \langle
A_s^{-1/2}w,A_s^{-1/2}w\rangle^{1/2} \\&=\langle
A^{-1}v,A_sA^{-1}v\rangle^{1/2}\langle A_s^{-1}w,w\rangle^{1/2}\\ &
=\langle A^{-1}v,v\rangle^{1/2}\langle
A_s^{-1}w,w\rangle^{1/2}\\&=\langle (A^{-1})_s
v,v\rangle^{1/2}\langle A_s^{-1}w,w\rangle^{1/2}\\& \leq\langle
A_s^{-1}v,v\rangle^{1/2}\langle A_s^{-1}w,w\rangle^{1/2}.\\
\end{split} \end{eqnarray} Then
\begin{equation}\label{eq2.2} \langle\Gamma v,w\rangle=\langle A^{-1}Bv,Bw\rangle\leq\langle A_s^{-1}Bv,Bv\rangle^{1/2}\langle A_s^{-1}Bw,Bw\rangle^{1/2} .\end{equation}\\
By Cauchy-Schwarz inequality, we have
\begin{equation}\label{eq2.3}
\langle Dv,w\rangle\leq\langle Dv,v\rangle^{1/2}\langle
Dw,w\rangle^{1/2}.
\end{equation}
Combining \eqref{eq2.2} and \eqref{eq2.3} and using Cauchy-Schwarz
inequality again we get
\begin{equation}\label{eq2.40}
\begin{split}
\langle(\Gamma+D)v,w\rangle&\leq\langle(B^TA_s^{-1}B+D)v,v\rangle^{1/2}\langle(B^TA_s^{-1}B+D)w,w\rangle^{1/2}\\&=\Vert
v\Vert_H\Vert w\Vert_H.
\end{split}
\end{equation}
According to Lemma 2.2, we obtain
\begin{equation*}
\frac{\langle Hv,v\rangle}{\langle G_{i}v,v\rangle}=\frac{\langle
G^{-1}_{i}Hv,v\rangle}{\langle v,v\rangle}\leq \theta(1+\beta_1),
\end{equation*}
i.e.,\\
\begin{equation*}
\langle Hv,v\rangle\leq \theta(1+\beta_1)\langle G_{i}v,v\rangle,
\end{equation*}
or equivalently
\begin{equation}\label{eq2.41}
\Vert v \Vert_H^2\leq\theta(1+\beta_1)\Vert v \Vert_{G_{i}}^2.
\end{equation}
Combining \eqref{eq2.40} and  \eqref{eq2.41}, we get
\begin{equation}\label{eq2.4}
\langle(\Gamma+D)v,w\rangle\leq \theta(1+\beta_1)\Vert v
\Vert_{G_{i}}\Vert w \Vert_{G_{i}}.
\end{equation}
Let $ w=G_{i}^{-1}(\Gamma+D)v $ and substitute $w$ into
\eqref{eq2.4}, we obtain
\begin{equation}\label{eq2.5}
\langle(\Gamma+D)v,{G^{-1}_{i}}(\Gamma+D)v)\rangle\leq(\theta(1+\beta_1))^2\Vert
v\Vert_{G_{i}}^2.
\end{equation}
 On the other hand, according to Lemma 2.1, Lemma 2.2 and $\alpha\geq1$ gives \begin{eqnarray}\label{eq2.6}\langle(\Gamma+D) u,u\rangle\geq\frac{1}{\alpha^2}(\langle A_s^{-1}Bu,Bu\rangle+\langle Du,u\rangle)\geq \frac{\theta(1-\beta_1)}{\alpha^2}\Vert u\Vert_{G_{i}}^2.
\end{eqnarray}
According to \eqref{eq2.1}, \eqref{eq2.5} and \eqref{eq2.6}, we have
\begin{equation*}
 ||(I-\theta\tau_i\widehat{S}^{-1}(B^{T}A^{-1}B+D))u||^{2}_{G_{i}}
 \leq ( 1-\frac{2\theta(1-\beta_1)}{\alpha^2}+\theta^2(1+\beta_1)^2 ) ||u||^2_{G_{i}}.\\
\end{equation*}\\
If we choose
\begin{equation*}
0<\theta<\frac{1-\beta_1}{\alpha^2(1+\beta_1)^2},
\end{equation*}
then by simple manipulations we have
\begin{equation*}
 ||(I-\theta\tau_i\widehat{S}^{-1}(B ^{T}A^{-1}B+D))u||^{2}_{G_{i}}\leq(1-\frac{\theta(1-\beta_1)}{\alpha^2})||u||^2_{G_{i}}.
 \end{equation*}
This concludes the proof of the theorem. ~ ~ ~  ~~~~~~~~~~~~~~~~~~~~~~~~~~~~~~~~~~~~~~~~~~~~~~~~~~~~~~~~~~~~~~~~~~~~~~~~~~~~~~~~~~~~~$\Box$ \\
{\bf Remark}: Some remarks about the choice of $\tau_i$. \\
\indent 1. If we substitute $H$ in \eqref{eq2.01} by $S_s$, we can
suggest another strategy to choose variable parameter $\tau_i$. That
is, the relaxation parameter $\tau_i$ can be determined such that
the norm
\begin{equation}
\Vert\tau_i\widehat S^{-1}g_i- S_s^{-1}g_i\Vert_{S_s}^2
 \end{equation}
 is minimized, where $g_i=B^Tx_{i+1}-Dy_i-g$, and $S_s$ is the spd part of $S$, then\\
 \begin{equation}\label{eq2.12}
 \tau_i= \left\{
   \begin{array}{cc}
  \frac{\langle g_i,\widehat {S}^{-1}g_i\rangle}{\langle(B^T(A^{-1})_sB+D)\widehat {S}^{-1}g_i,\widehat {S}^{-1}g_i\rangle},&g_i\neq0  ,\\
   1,&g_i=0 .\\
   \end{array}
   \right.
  \end{equation}
\indent 2. Similar to Lemma 2.2, we can derive the following
conclusion.

For any natural number $i$, there is a symmetric and positive definite ${m \times m}$ matrix ${ G_{i} }$ such that\\
(i) $G_{i}^{-1}g_i=\theta\tau_i\widehat {S}^{-1}g_i $ with $g_i=B^Tx_{i+1}-Dy_i-g $ as defined in Algorithm 2.1;\\
(ii) All eigenvalues of the matrix $ G_{i}^{-1}S_s $ lie in the interval $[\theta(1- \beta_2),\theta(1+\beta_2)]$, where $\beta_2=\frac{\kappa_2-1}{\kappa_2+1}, \kappa_2=cond(\widehat{S}^{-1}S_s), S_s=B^T(A^{-1})_sB+D$. \\
\indent 3. With the statements above, we can derive the similar
convergence result for the parameter choice strategy \eqref{eq2.12}.
We just need to apply Lemma 2.1 again on \eqref{eq2.02}, that is to
say, we can obtain
 \begin{equation*}
 \langle A^{-1}v,w\rangle \leq\langle A_s^{-1}v,v\rangle^{1/2}\langle A_s^{-1}w,w\rangle^{1/2}\leq \alpha^2\langle (A^{-1})_sv ,v\rangle^{1/2}\langle (A^{-1})_sw ,w\rangle^{1/2}.
 \end{equation*}
 Then using the same strategy as the proof of Theorem 2.1, we can obtain that, when we choose $0<\theta<\frac{1-\beta_2}{\alpha^4(1+\beta_2)^2}$, we have\\
 \begin{equation*}
 ||(I-\theta\tau_i\widehat {S}^{-1}(B^TA^{-1}B+D))u||^{2}_{G_{i}}\leq(1-\theta(1-\beta_2))||u||^2_{G_{i}}, ~\forall \ u\in{H_2}.
  \end{equation*}

 \section{Analysis of  linear inexact Uzawa algorithm with relaxation parameter}

\indent In this section we define and analyze a linear one-step
Uzawa type algorithm with relaxation parameter for \eqref{eq1.0}.
Under the minimal assumption needed to guarantee solvability, we
suggest an efficient and simple method for solving \eqref{eq1.0} .
The exact inverse of $A$ is replaced by a preconditioner $A_0$ for
the symmetric part of $A$. Let $A_0$ be a linear, symmetric positive
definite operator and satisfy \eqref{eq1.11}.

 {\bf Algorithm 3.1} Given $x_0\in
{H_1}$ and $y_0\in {H_2}$, the sequence$\{(x_i,y_i)\}$ is defined,
for ${i=0,1,2,...,}$ by
\begin{eqnarray*} x_{i+1} & = & x_i+\omega A_0^{-1}(f-(Ax_i+By_i)),
\\
y_{i+1}&=&y_i+\delta\tau_i\widehat{S}^{-1}(B^Tx_{i+1}-Dy_i-g),
\end{eqnarray*}
where \begin{equation} \label{eq3.00} \tau_i= \left\{
   \begin{array}{cc}
  \frac{\langle g_i,\widehat{S}^{-1}g_i\rangle}{\langle(B^TA_0^{-1}B+D)\widehat{S}^{-1}g_i,\widehat{S}^{-1}g_i\rangle},&g_i\neq0,  \\
   1,&g_i=0. \\
   \end{array}
   \right.\end{equation}
 Here $\omega$ and $\delta$ are positive constant parameters determined to guarantee the convergence, $\tau_i$ above can be computed effectively as the method in \cite{HuZou_MAA01}, while we work on the nonsymmetric matrix $A$ and $D\neq 0$. We will assume that $\omega<1/\kappa_0$. It then follows from \eqref{eq1.11} that $A_0-\omega A_s$ is positive definite.

{\bf Theorem 3.1} Suppose that $A$ has a positive definite symmetric part $A_s$, satisfying \eqref{eq1.1}. Suppose also that $A_0$ is symmetric positive definite operator satisfying \eqref{eq1.11}. Then Algorithm 3.1 is convergent if $\delta<1/2$, $0<\omega<min(\frac{1}{3\alpha^2\kappa_0^2},\frac{1+\kappa_0(1-\delta(1+\beta_3))}{(\alpha^2\kappa_0+1)\kappa_0})$. Moreover, when $\delta<\frac{1}{4\alpha^2\kappa_0^2}$, the iteration errors $e_i^x$ and $e_i^y$ satisfying\\
\begin{equation}
\left(\omega^{-1}\Vert e_i^x\Vert_{A_0-\omega A_s}^2+\Vert
e_i^y\Vert_{G_{i}}^2\right)^{1/2}\leq\bar\rho^{i}\left(\omega^{-1}\Vert
e_0^x\Vert_{A_0-\omega A_s}^2+\Vert
e_0^y\Vert_{G_{0}}^2\right)^{1/2},
\end{equation}
for any $i\geq1$. Here
\begin{equation}
\bar\rho=\frac{\omega/2-\omega\Delta+\sqrt{(\omega/2-\omega\Delta)^2+4(1-\omega/2)}}{2},
\end{equation}
where $\Delta=\frac{\delta(1-\beta_3)}{\kappa_0}$.\\

{\bf Lemma 3.1} With  the assumption of \eqref{eq1.11}, for any natural number $i$, there is a symmetric and positive definite ${m \times m}$ matrix ${G_{i}}$  such that\\
(i) $G_{i}^{-1}g_i=\delta\tau_i\widehat{S}^{-1}g_i $ with $g_i=B^Tx_{i+1}-Dy_i-g $ as defined in Algorithm 3.1;\\
 (ii) All eigenvalues of the matrix $ G_{i}^{-1}H$ lie in the interval $[\frac{\delta(1-\beta_3)}{\kappa_0},\delta(1+\beta_3)]$, where $\beta_3=\frac{\kappa_3-1}{\kappa_3+1}, \kappa_3=cond(\widehat{S}^{-1}M), M=B^TA_0^{-1}B+D$. \\
The proof of this lemma is similar to  Lemma 3.2 in
\cite{HuZou_MAA01}.

   In order to analyze Algorithm 3.1 we formulated it in terms of the iteration errors. It is easy to see that $e_i^{x}$ and $e_i^y$ satisfy the following equations.
   \begin{eqnarray*}
   \begin{split}
   e_{i+1}^x&=e_i^x-\omega A_0^{-1} (Ae_i^x+Be_i^y),\\
   e_{i+1}^y&=(I-\omega\delta\tau_i\widehat{S}^{-1}(B^TA_0^{-1}B+D))e_i^y+\delta\tau_i\widehat{S}^{-1}B^T(I-\omega A_0^{-1}A)e_i^x.
  \end{split}
   \end{eqnarray*}
   From Lemma 3.1, we have
   \begin{eqnarray*}
   \begin{split}
   e_{i+1}^x&=e_i^x-\omega A_0^{-1} (Ae_i^x+Be_i^y),\\
   e_{i+1}^y&=(I-\omega G_{i}^{-1}(B^TA_0^{-1}B+D))e_i^y+G_{i}^{-1}B^T(I-\omega A_0^{-1}A)e_i^x.
  \end{split}
   \end{eqnarray*}
   For convenience, these equations can be written in the matrix form as\\
$$\begin{pmatrix}
  e_{i+1}^x \\e_{i+1}^y   \end{pmatrix}
= \begin{pmatrix} I-\omega A_0^{-1}A & -\omega
A_0^{-1}B\\G_{i}^{-1}B^{T}(I-\omega A_0^{-1}A) & I-\omega
G_{i}^{-1}(B^{T}A_0^{-1}B+D) \end{pmatrix}
\begin{pmatrix}
  e_{i}^x \\e_{i}^y
\end{pmatrix}. $$
Then straightforward manipulation yields
\begin{equation*}
\mathcal{N}e_{i+1}=\mathcal{M}e_i ,
\end{equation*}
where
\begin{equation*}
e_i= \begin{bmatrix} e_i^x \\
 e_i^y \end{bmatrix},
 \end{equation*}
 \begin{equation*}
 \mathcal{N}=\begin{bmatrix}
 \omega^{-1}(A_0-\omega A^T) &0\\0& {G_{i}}
 \end{bmatrix},
  \end{equation*}
\begin{equation*}
 \mathcal{M}=\begin{bmatrix}
 \omega^{-1}(A_0-\omega A^T)A_0^{-1}(A_0-\omega A) & -(A_0-\omega A^T)A_0^{-1}B\\B^TA_0^{-1}(A_0-\omega A)&{G_{i}}-\omega (B^TA_0^{-1}B+D)
 \end{bmatrix}.
  \end{equation*}
It is clear that we can study the convergence of Algorithm 3.1 by investigating the properties of the linear operator $ \mathcal{M} $ and $\mathcal{N}$. We shall reduce this problem to estimate of the spectral radius of related symmetric operators.\\
Let $ \mathcal{M}_1$ be the symmetric matrix defined by
\begin{equation*}
\mathcal{M}_1=\mathcal{J}\mathcal{M},
\end{equation*}
where
\begin{equation*}
\mathcal{J}= \left(
   \begin{array}{cc}
   -I& 0\\0& I
   \end{array}
\right).
\end{equation*}

Let $\mathcal{N}_s$ be the symmetric part of $\mathcal{N}$. Since
$\omega < 1/\beta$, $\mathcal{N}_s$ is symmetric positive definite.
The proof of the convergence
 needs the   estimation of the eigenvalues of the following
generalized eigenvalue problem
\begin{equation} \label{DefineGeneralizedEig4ConvergenceProof}
\lambda \mathcal{N}_s \psi = \mathcal{M}_1 \psi.
\end{equation}
%Since $\omega < 1/\beta$, $\mathcal{N}_s$ is symmetric positive definite and the above problem is well defined.
%Because $\mathcal{N}_s$ and $M_1$ are symmetric, the eigenvalues $\lambda$ are real.
%
Let $\{ (\lambda_i, \psi_i) \} $ be the eigenpairs for
\eqref{DefineGeneralizedEig4ConvergenceProof}  and  $\langle
\mathcal{N}_s \psi_i, \psi_j \rangle = \delta_{ij}$.
 Any vectors $v$ and $w$ in $H_1 \times H_2$ can be represented as $v =
\Sigma_i v_i \psi_i$ and $w = \Sigma_j w_j \psi_j$, and hence
\cite{BramblePasciakVassilev_MC99},
\begin{eqnarray*}
\langle \mathcal{M}_1 v, w \rangle = \Sigma_{ij} v_i w_j \langle
\mathcal{M}_1 \psi_i,
\psi_j \rangle = \Sigma_j v_j w_j \lambda_j \\
\leq \bar{\rho} \sqrt{\Sigma_j v_j^2} \sqrt{\Sigma_j w_j^2} =
\bar{\rho} ||v||_{\mathcal{N}_s} ||w||_{\mathcal{N}_s}.
\end{eqnarray*}
It is easy to verify that \cite{BramblePasciakVassilev_MC99}
\begin{eqnarray*}
\langle \mathcal{N}_s e_{i+1}, e_{i+1} \rangle & = & \langle
\mathcal{M} e_i, e_{i+1} \rangle = \langle \mathcal{J} \mathcal{M}_1
e_i, e_{i+1} \rangle =  \langle \mathcal{M}_1 e_i,
\mathcal{J}e_{i+1} \rangle \\
&\leq& \bar{\rho} || e_i||_{\mathcal{N}_s}
||\mathcal{J}e_{i+1}||_{\mathcal{N}_s} = \bar{\rho} ||
e_i||_{\mathcal{N}_s} || e_{i+1}||_{\mathcal{N}_s}.
\end{eqnarray*}

In order to prove Theorem 3.1, we need the following two lemma.

{\bf Lemma 3.2} The iteration error $e_i$ satisfies
\begin{equation}
\langle \mathcal{N}_se_{i+1},e_{i+1}\rangle^{1/2}\leq\bar\rho\langle
\mathcal{N}_se_i,e_i\rangle^{1/2}.
\end{equation}
where  $\bar\rho= \max\vert \lambda_i\vert $, with\{$\lambda_i\}$
the eigenvalues of \eqref{DefineGeneralizedEig4ConvergenceProof}.
See Lemma 3.1 in \cite{BramblePasciakVassilev_MC99}.

%\indent Then from Lemma 3.2, we know that in order to prove the convergence of the Theorem 3.1, we just need to estimate the eigenvalues of the following generalized eigenvalue problem\\
%\begin{equation}\label{eq3.1}
%\lambda N_s\varphi =M_1\varphi,
%\end{equation}
%since $N_s$ and $M_1$ are symmetric, the eigenvalues $\lambda$ are real.

 {\bf Lemma 3.3} Let $A_0$  satisfy \eqref{eq1.11} and $\omega$
be a positive number with $\omega<\frac{1}{\kappa_0}$ ,
then \\
\begin{equation*}
\Vert(I-\omega A_0^{-1}A)v\Vert_{A_0}^2\leq\bar\omega\langle(A_0-\omega A_s)v,v\rangle,\\
\end{equation*}
where\\
\begin{equation*}
\bar\omega=1-\omega+\frac
{\omega^2\alpha^2\kappa_0^2}{1-\omega\kappa_0}.
\end{equation*}
For the proof of this lemma please see Lemma 3.2 in
\cite{BramblePasciakVassilev_MC99}. \\

\emph{Next we prove the convergence of Algorithm 3.1}. The proof is
analogous to that of Theorem 3.1 in
\cite{BramblePasciakVassilev_MC99}. Because of Lemma 3.2, it
suffices to bound the eigenvalue of
\eqref{DefineGeneralizedEig4ConvergenceProof}. We begin with the
negative eigenvalues. Let $(\chi,\xi)$ be an eigenvector with
eigenvalue $\lambda<0$ . Then multiplying the first equation of
\eqref{DefineGeneralizedEig4ConvergenceProof} by $ A_0(A_0-\omega
A^T)^{-1}$ gives
\begin{equation}\label{eq3.2}
\lambda\omega^{-1}A_0(A_0-\omega A^T)^{-1}(A_0-\omega
A_s)\chi=-\omega^{-1}(A_0-\omega A)\chi+B\xi.
\end{equation}
The second equation of \eqref{DefineGeneralizedEig4ConvergenceProof}
is
\begin{equation}\label{eq3.3}
\lambda G_{i}\xi=B^TA_0^{-1}(A_0-\omega A)\chi+({G_{i}}-\omega
(B^TA_0^{-1}B+D))\xi.
\end{equation}
Applying $\omega B^TA_0^{-1} $ to \eqref{eq3.2} and adding it to
\eqref{eq3.3}, we obtain
\begin{equation}\label{eq3.4}
((1-\lambda) {G_{i}}-\omega D)\xi=\lambda B^T(A_0-\omega
A^T)^{-1}(A_0-\omega A_s)\chi.
\end{equation}
We note that $\lambda<0$, then
\begin{equation}\label{eq3.5}
\langle((1-\lambda) {G_{i}}-\omega
D)v,v\rangle\geq\langle((1-\lambda) {G_{i}}-
D)v,v\rangle\geq(1-\lambda)\langle( {G_{i}}- D)v,v\rangle.
\end{equation}
From lemma 3.1,
\begin{equation}
\langle(B^TA_s^{-1}B+D)v,v\rangle\leq\delta(1+\beta_3)\langle
{G_{i}}v,v\rangle\leq\langle {G_{i}}v,v\rangle.
\end{equation}
Thus
\begin{equation*}
\langle( {G_{i}}- D)v,v\rangle\geq\langle B^TA_s^{-1}Bv,v\rangle>0.
\end{equation*}
That is to say, \eqref{eq3.5} implies $(1-\lambda) {G_{i}}-\omega D$
is positive. Then we can eliminate $\xi$ among \eqref{eq3.2} and
\eqref{eq3.4}, we find that
\begin{equation*}
\begin{split}
-\frac{1}{\lambda}(A_0-\omega A)\chi+&\omega B((1-\lambda)
{G_{i}}-\omega D)^{-1}B^T(A_0-\omega A^T)^{-1}(A_0-\omega
A_s)\chi\\&=A_0(A_0-\omega A^T)^{-1}(A_0-\omega A_s)\chi.
\end{split}
\end{equation*}
Taking the inner product with $(A_0-\omega A^T)^{-1}(A_0-\omega
A_s)\chi$ yields
\begin{equation}
\begin{split}
-\frac{1}{\lambda}((A_0-\omega A_s)\chi,\chi) +& \omega\Vert
B^{T}(A_0-\omega A^{T})^{-1}(A_0-\omega A_s)\chi\Vert_{((1-\lambda)
{G_{i}}-\omega D)^{-1}}^2\\ =& \Vert (A_0-\omega
A^{T})^{-1}(A_0-\omega A_s)\chi\Vert_{A_0}^2.
\end{split}
\end{equation}
For convenience, the above equation can be abbreviated as\\
\begin{equation*}
\Gamma_1+\Gamma_2=\Gamma_3.
\end{equation*}
In order to bound $\Gamma_2$ we note that for any $\phi\in {H_1}$,\\
\begin{equation*}
\begin{split}
\langle((1-\lambda) {G_{i}}-\omega D)^{-1}B^T\phi,B^T\phi\rangle =& \sup_{\zeta\in H_2}\frac{\langle\phi,B\zeta\rangle^2}{\langle((1-\lambda) {G_{i}}-\omega D)\zeta,\zeta\rangle}\\
=& \sup_{\zeta\in H_2}\frac{\langle(A_s)^{1/2}\phi,(A_s)^{-1/2}B\zeta\rangle^2}{\langle((1-\lambda){G_{i}}-\omega D)\zeta,\zeta\rangle}\\
\leq & \sup_{\zeta\in H_2}\frac{\langle A_s\phi,\phi\rangle\langle
B^{T}A_s^{-1}B\zeta,\zeta\rangle}{\langle((1-\lambda) {G_{i}}-\omega
D)\zeta,\zeta\rangle}\\ \leq & \frac{1}{1-\lambda}\frac{\langle
A_s\phi,\phi\rangle\langle
B^{T}A_s^{-1}B\zeta,\zeta\rangle}{\langle({G_{i}}-D)\zeta,\zeta\rangle}
\\ \leq &\frac{1}{1-\lambda}\delta(1+\beta_3)\langle
A_s\phi,\phi\rangle.\\
\leq &\frac{\delta(1+\beta_3)\kappa_0}{1-\lambda}\langle
A_0\phi,\phi\rangle.
\end{split}
\end{equation*}
Therefore,\\
\begin{equation*}
\Gamma_2\leq\omega\frac{\delta(1+\beta_3)\kappa_0}{1-\lambda}\Gamma_3=\frac{\omega\delta(1+\beta_3)\kappa_0}{(1-\lambda)}\Gamma_3.
\end{equation*}
By Lemma 3.3, we have\\
\begin{equation*}
\begin{split}
((A_0-\omega A^{T})A_0^{-1}(A_0-\omega A)\phi,\phi)&=\Vert(I-\omega
A_0^{-1}A)\phi\Vert_{A_0}^2\\ &\leq\bar\omega\langle(A_0-\omega
A_s)\phi,\phi\rangle,
\end{split}
\end{equation*}
which is equivalent to\\
\begin{equation*}
\langle(A_0-\omega A_s)^{-1}\phi,\phi\rangle\leq\bar\omega
\langle(A_0-\omega A)^{-1}A_0(A_0-\omega
A^{T})^{-1}\phi,\phi\rangle.
\end{equation*}
Taking $ \phi=\langle(A_0-\omega A_s)\chi,\chi\rangle$ in the formula above, we obtain\\
\begin{equation*}
\Gamma_3\geq\frac{1}{\bar\omega}\langle(A_0-\omega
A_s)\chi,\chi\rangle.
\end{equation*}
Then using the fact that $\lambda \leq0$, we have
\begin{equation}
\begin{split}
-\frac{1}{\lambda}\langle(A_0-\omega
A_s)\chi,\chi\rangle&\geq\frac{1}{\bar\omega}(1-\frac{\omega\delta(1+\beta_3)\kappa_0}{1-\lambda})\langle(A_0-\omega
A_s)\chi,\chi\rangle\\& \geq
\frac{1-\omega\delta(1+\beta_3)\kappa_0}{\bar\omega}\langle(A_0-\omega
A_s)\chi,\chi\rangle.
\end{split}
\end{equation}
Thus\\
\begin{equation}
-\lambda \leq \frac{\bar\omega}{1-\omega\delta(1+\beta_3)\kappa_0}.
\end{equation}
Indeed, by simple manipulations, when
\begin{equation}
\omega<
\frac{1+\kappa_0(1-\delta(1+\beta_3))}{(\alpha^2\kappa_0+1)\kappa_0},
\end{equation}
we have
\begin{equation*}
-\lambda <1.
\end{equation*}
In addition, when\begin{equation*}
0<\omega<\frac{1}{3\alpha^2\kappa_0^2},
\end{equation*}
we get
\begin{equation} \label{eqn4baromega}
\bar\omega=1-\omega+\frac{\alpha^2\kappa_0^2\omega^2}{1-\omega\kappa_0}\leq1-\omega(1-\frac{1/3}{1-1/3})=1-\omega/2.
\end{equation}
When $\delta<\frac{1}{4(1+\beta_3)}$, then
\begin{equation*}
1-\omega\delta(1+\beta_3)\geq1-\omega/4,
\end{equation*}
thus \\
\begin{equation*}
-\lambda\leq\frac{1-\omega/2}{1-\omega/4}<1-\omega/4,
\end{equation*}
which provides a bound for the negative part of the spectrum.\\

Next we obtain a bound for the positive eigenvalues.
 We set
$S_0=B^TA_0^{-1}B+D$ , to this end we factor $\mathcal{M}_1$ as
\begin{equation*}
\mathcal{M}_1=\mathcal{P}^{T}\mathcal{M}_2\mathcal{P},
\end{equation*}
where
\begin{equation*}
\mathcal{P}=\left(
\begin{array}{cc}
\theta^{-1/2}A_0^{-1/2}(A_0-\omega A) & 0\\0 & I
  \end{array} \right),
\end{equation*}
\begin{equation*}
\mathcal{M}_2=\left(
\begin{array}{cc}
-\omega^{-1}\theta I & \theta^{1/2}A_0^{-1/2}B\\
\theta^{1/2}B^{T}A_0^{-1/2}& G_{i}-\omega S_0
 \end{array}
 \right),
\end{equation*}
and $\theta>0$. The largest eigenvalue is given by\\
\begin{equation*}
\Lambda= \sup_{w\in{H_1\times H_2}} \frac{\langle
\mathcal{M}_1w,w\rangle}{\langle \mathcal{N}_sw,w\rangle}=
\sup_{w\in{H_1\times H_2}} \frac{\langle
\mathcal{M}_2\mathcal{P}w,\mathcal{P}w\rangle}{\langle
\mathcal{N}_sw,w\rangle}.
\end{equation*}
In order to obtain an upper bound for $\Lambda$, from Lemma 3.3 we
note that
\begin{equation*}
\begin{split}
\Vert(I-\omega A_0^{-1}A)\chi\Vert_{A_0}^2=\Vert
A_0^{-1/2}(A_0-\omega A)\chi \Vert^2 &\leq \bar\omega
\langle(A_0-\omega A_s)\chi,\chi\rangle.
\end{split}
\end{equation*}
We especially choose $\theta=1-\omega/2$ and use
\eqref{eqn4baromega}, then we have
$$
\Vert(I-\omega A_0^{-1}A)\chi\Vert_{A_0}^2
\leq\theta\langle(A_0-\omega A_s)\chi,\chi\rangle. $$ That is,
\begin{equation*}
\theta^{-1}\Vert A_0^{-1/2}(A_0-\omega A)\chi
\Vert^2\leq\langle(A_0-\omega A_s)\chi,\chi\rangle.
\end{equation*}
And thus
\begin{equation*}
\langle
\begin{bmatrix}\omega^{-1}I&0 \\0& G_i \end{bmatrix}
\mathcal{P}\begin{bmatrix} \chi \\ \xi   \end{bmatrix} ,
\mathcal{P}\begin{bmatrix} \chi
\\ \xi \end{bmatrix}\rangle = \theta^{-1}\omega^{-1}\Vert
A_0^{-1/2}(A_0-\omega A)\chi \Vert^2+\Vert \xi\Vert_{G_{i}}^2 \leq
\langle \mathcal{N}_s\begin{bmatrix}\chi
\\ \xi \end{bmatrix},\begin{bmatrix}\chi \\ \xi  \end{bmatrix}\rangle.
\end{equation*}
Thus, it suffices to show that for any vector
$(\phi,\zeta)\in{H_1\times H_2}$,
\begin{equation*}
\begin{split}
\langle \mathcal{M}_2\mathcal{P}\begin{bmatrix}\phi \\ \zeta
\end{bmatrix} ,\mathcal{P}\begin{bmatrix}\phi \\ \zeta
\end{bmatrix}\rangle&\leq\bar \rho \ \langle
\begin{bmatrix}\omega^{-1}I&0 \\0& G_i \end{bmatrix}
\mathcal{P}\begin{bmatrix}\phi \\ \zeta  \end{bmatrix} ,
\mathcal{P}\begin{bmatrix}\phi
\\ \zeta  \end{bmatrix}\rangle.
\end{split}
\end{equation*}
Then $\bar\rho$ will be an upper bound for $\Lambda$. To this end
let $L=A_0^{-1/2}B$. Now $\mathcal{M}_2$ can be written as
\begin{equation*}
\mathcal{M}_2=\left(
\begin{array}{cc}

-\omega^{-1}\theta I & \theta^{1/2}L\\ \theta^{1/2}L^{T}&
G_{i}-\omega( L^{T}L+D)

 \end{array}
 \right).
\end{equation*}
Then we reduced the problem to estimate the largest eigenvalue
$\lambda$, with eigenvector $(\chi,\xi)$ satisfying
\begin{equation}\label{eq3.1.0}
\lambda  \begin{bmatrix}\omega^{-1}I&0 \\0& G_{i} \end{bmatrix}
\begin{bmatrix}\chi \\ \xi
\end{bmatrix}=\mathcal{M}_2\begin{bmatrix}\chi \\ \xi
\end{bmatrix}.
\end{equation}
The first equation of \eqref{eq3.1.0} is
\begin{equation}\label{eq3.1.1}
-\theta\omega^{-1}\chi+\theta^{1/2}L\xi=\lambda\omega^{-1}\chi,
\end{equation}
and the second is
\begin{equation}\label{eq3.1.2}
\theta^{1/2}L^{T}\chi+(G_{i}-\omega L^{T}L-\omega D)\xi=\lambda
G_{i}\xi.
\end{equation}
Solving for $\chi$ in \eqref{eq3.1.1}, we get
\begin{equation*}
\chi=\omega(\lambda+\theta)^{-1}\theta^{1/2}L\xi.
\end{equation*}
Substituting this in \eqref{eq3.1.2} yields
\begin{equation}\label{eq3.1.3}
\begin{split}
(1-\lambda)(\lambda+\theta)\langle
G_{i}\xi,\xi\rangle&=\omega\lambda \langle
L\xi,L\xi\rangle+(\lambda+\theta)\omega\langle
D\xi,\xi\rangle\\&\geq\omega\lambda\langle(L^TL+D)\xi,\xi\rangle.
\end{split}
\end{equation}
In addition,
\begin{equation}\label{eq3.1.4}
\begin{split}
\langle(L^{T}L+D)\xi,\xi\rangle&=\langle
A_0^{-1}B\xi,B\xi\rangle+\langle D\xi,\xi\rangle\\&\geq\langle
A_s^{-1}B\xi,B\xi\rangle+\langle
D\xi,\xi\rangle\\&\geq\frac{\delta(1-\beta_3)}{\kappa_0}\langle
G_{i}\xi,\xi\rangle,
\end{split}
\end{equation}
where the last inequality is derived by using Lemma 3.1.\\
Then from \eqref{eq3.1.3} and \eqref{eq3.1.4}, we have
\begin{equation*}
(1-\lambda)(\lambda+\theta)\geq\omega\lambda\frac{\delta(1-\beta_3)}{\kappa_0},
\end{equation*}
or equivalently
\begin{equation}
\lambda^2-\lambda(1-\theta-\omega\frac{\delta(1-\beta_3)}{\kappa_0})-\theta\leq0.
\end{equation}
From here we obtain that
\begin{equation}
\begin{split}
\lambda&\leq\frac{1-\theta-\omega\Delta+\sqrt{(1-\theta-\omega\Delta)^2+4\theta}}{2}\\&=\frac{\omega/2-\omega\Delta+\sqrt{(\omega/2-\omega\Delta)^2+4(1-\omega/2)}}{2},
\end{split}
\end{equation}
where $\Delta=\frac{\delta(1-\beta_3)}{\kappa_0}$.\\

When $\delta<1/2, \frac{\delta(1-\beta_3)}{\kappa_0}\leq 1/4$, using
the simple algebraic manipulation, similar to Remark 3.2 in
\cite{BramblePasciakVassilev_MC99}, it follows that
\begin{equation}
\frac{\omega/2-\omega\Delta+\sqrt{(\omega/2-\omega\Delta)^2+4(1-\omega/2)}}{2}\leq1-\frac{\omega\delta(1-\beta_3)}{2\kappa_0}<1,
\end{equation}
which provides a bound for the positive part of the spectrum. \\

Finally, elementary inequalities \cite{BramblePasciakVassilev_MC99}
imply that
\begin{equation*}
1-\frac{\omega}{4}\leq\frac{\omega/2-\omega\Delta+\sqrt{(\omega/2-\omega\Delta)^2+4(1-\omega/2)}}{2},
\end{equation*}
which concludes the proof of the theorem.~~~~~~~~~~~~~~~~~~~~~~~~~~~~~~~~~~~~~~~~~~~~~~~~~~~~~~~~~~~~~~~~~~~~~~~~~~~~~~~~~~~~$\Box$\\

{\bf Remark:} The convergence of Algorithm 3.1 for the saddle point
problems \eqref{eq1.0} under the condition \eqref{eq1.11}, i.e., the
preconditioner $A_0$ for $A_s$ is appropriately scaled. In \cite{K.
H. J_MAA12 }, it shows that the assumption \eqref{eq1.11} can be
removed in the convergence proof of Algorithm 1.1. This can be
applied on our Algorithm 3.1 also, that is, the convergence of
Algorithm 3.1 without the condition \eqref{eq1.11} can be achieved.

\section{Numerical examples }
In this section, we present some numerical experiments to show the
performance of Algorithm 3.1 with parameter $\omega=0.3$ ,
$\delta=0.3$ and $\tau_i$ selected by \eqref{eq3.00}. The numerical
tests are performed by \textsc{Matlab} R2008a on the laptop with
Intel(R) Core(TM) i5-3210M CPU@ 2.50GHZ and 4G RAM. We first
consider Oseen problem  in the rectangular domain
$\Omega=(0,1)\times(0,1)$, with Dirichlet boundary conditions:
$u=v=0$ on $x=0, x=1, y=0; u=1, v=0,$ on $y=1$. The governing
equations
 are
 \begin{eqnarray*}
 -\nu\Delta u +(w \cdot\nabla)+\nabla p&=&f, \\
-\nu\Delta u &=&0,
\end{eqnarray*}
¡¡where $w=(w_1,w_2)$ denote the wind. In this paper, we choose $w$,
such that it is the image of $[2y(1-x^2),-2x(1-y^2)]$ under the
mapping from $(-1,1)\times (-1,1)$ to $\Omega$. That is,
\begin{equation*}
 w=
   \begin{bmatrix}
  8x(1-x)(2y-1) \\
   -8y(1-y)(2x-1) \\
   \end{bmatrix}.
  \end{equation*}

\indent  We take different preconditioners to test the convergence of  the Algorithm 3.1. The preconditioner for $A_0$ is chosen from one of the following four: the \textsc{Jacobi} preconditioner, the \textsc{Ilu(1e-1)} preconditioner, the incomplete Cholesky factorization with the drop tolerance of $10^{-1}$, and the exact symmetric part of (2,2) block, i.e, $A_s$. We consider the case where $\nu=1$ with $32\times32$ grids. And the preconditioner $\widehat S$ for Schur complement is approximated by a scaled identity matrix.\\
\begin{figure}[!htbp]
  \centering
\includegraphics [width=10cm,height=4.5cm]{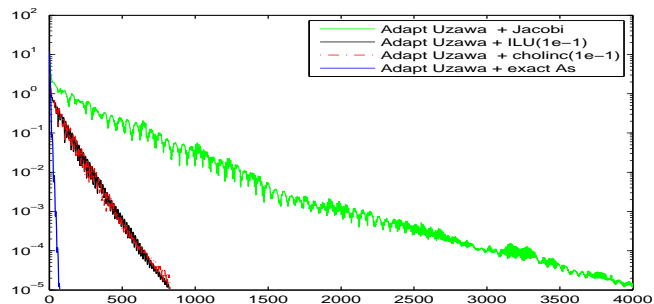}\caption{Residuals of adaptive Uzawa algorithm for Oseen problem with different preconditions ($32\times32$ grids,  $ \nu=1$) }
\label{fig2}
 \end{figure}
\indent  Figure 1 illustrated the residuals of Algorithm 3.1 for
Oseen problem with different preconditioners. Obviously, the
convergence of the residual curves confirms our analysis in Section
3, where we prove the convergence Theorem 3.1 for Algorithm 3.1. \\
\indent From this figure, we can observe that, according to the
number of iteration steps, the exact preconditioner $A_s$ is the
most fast. While the preconditioners \textsc{Cholinc} and
\textsc{Ilu} with tolerance of $10^{-1}$ take more iteration steps.
The preconditioners \textsc{Cholinc} and \textsc{Ilu} almost need
the same number of iterations, but the CPU time and flops in each
step of \textsc{Cholinc} preconditioner are less than those of the
preconditioner \textsc{Ilu}. The \textsc {Jacobi} preconditioner
takes even more iteration steps, but it sometimes needs less CPU
time, due to its simplicity, which can be seen from the test on N-S
equation in the following. Moreover, we observed that, when the
tolerance is between $[10^{-1}, 10^{-4}]$, the iteration numbers of
\textsc{Cholinc} and \textsc{Ilu} preconditioners are decreasing
with  smaller tolerance. When  the tolerance is less than or equal
to $10^{-4}$, the iteration numbers of preconditioners
\textsc{Cholinc} and \textsc{Ilu} are almost the same as the exact
preconditioner $A_s$.

\indent Next, we consider the Navier-Stokes problem. Solving the
Navier-Stokes problem corresponds to solve an Oseen problem in every
Picard iteration. Hence, the strategies designed for Oseen problem
can be applied to Navier-Stokes problem as well. We compare the
Algorithm 3.1  with the algorithm of Bramble-Pasciak-Vassilev (BPV)
with  $\delta=0.1$, $\tau=0.01$  in
\cite{BramblePasciakVassilev_MC99} and \textsc {Gmres} algorithm
with no preconditioning. Table I, Table II, Table III illustrated
the number of iterations and CPU times (seconds) of  three cases
with different preconditioners when  $\nu=0.01,0.1,$ and $ 1$,
 respectively.

\begin{table}[!htbp]
\begin{center}
\begin{tabular}{lcccccc}\hline
 & \multicolumn{2}{c}{$n$=16 } & & \multicolumn{2}{c}{ $n$=32} \\
 \cline{2-3} \cline{5-6}
Algorithms     & iter & CPU  &  & iter & CPU \\
\hline
\textsc{AdaptiveUzawa+Ilu(1e-4)} &  14 & 4.26 & &  16 &  57.76 \\
\textsc{AdaptiveUzawa+Cholinc(1e-4)} & 13 & 0.91 & & 15 & 21.00 \\
\textsc{AdaptiveUzawa+Jacobi} & 12 & 1.81 & & 12 & 73.50 \\
\textsc{AdaptiveUzawa+Exact} & 13 & 2.82 & & 13 &28.13\\
\textsc{BPV+Ilu(1e-4)} & 12 & 54.07 & & 12 & 343.40 \\
\textsc{BPV+Cholinc(1e-4)} & 12 & 6.99 & &  12 & 84.83\\
\textsc{BPV+Jacobi} & 12 & 21.39 & & 12 & 773.68 \\
\textsc{Gmres} & 11 & 14.17 & & 11 & 96.76 \\
\hline
\end{tabular}
\caption {Comparison of the computation time for N-S with
$\nu=0.01$}
\end{center}
\end{table}
\begin{table}[!htbp]
\begin{center}
\begin{tabular}{lcccccc}\hline
 & \multicolumn{2}{c}{$n$=16 } & & \multicolumn{2}{c}{ $n$=32} \\
 \cline{2-3} \cline{5-6}
Algorithms     & iter & CPU  &  & iter & CPU \\
\hline
\textsc{AdaptiveUzawa+Ilu(1e-4)} &  7 & 1.26 & &  14 &  34.59 \\
\textsc{AdaptiveUzawa+Cholinc(1e-4)} & 7 & 0.33 & & 14 & 14.96 \\
\textsc{AdaptiveUzawa+Jacobi} & 15 & 0.87 & & 12 & 32.68 \\
\textsc{AdaptiveUzawa+Exact} & 7 & 0.64 & & 14 &18.74\\
\textsc{BPV+Ilu(1e-4)} & 7 & 523.43 & & 14 & 153.94 \\
\textsc{BPV+Cholinc(1e-4)} & 7 & 3.05 & &  14 & 44.40\\
\textsc{BPV+Jacobi} & 7 & 10.94 & & 7 & 461.29 \\
\textsc{Gmres} & 6 & 2.16 & & 11 & 16.09 \\
\hline
\end{tabular}
\caption {Comparison of the computation time for N-S with $\nu=0.1$}
\end{center}
\end{table}
\begin{table}[!htbp]
\begin{center}
\begin{tabular}{lcccccc}\hline
 & \multicolumn{2}{c}{$n$=16 } & & \multicolumn{2}{c}{ $n$=32} \\
 \cline{2-3} \cline{5-6}
Algorithms     & iter & CPU  &  & iter & CPU \\
\hline
\textsc{AdaptiveUzawa+Ilu(1e-4)} &  20 & 2.13 & &  27 &  63.09 \\
\textsc{AdaptiveUzawa+Cholinc(1e-4)} & 16 & 0.70 & & 27 & 32.24 \\
\textsc{AdaptiveUzawa+Jacobi} & 192 & 5.39 & & 418 & 382.53 \\
\textsc{AdaptiveUzawa+Exact} & 16 & 1.12 & & 28 &40.65\\
\textsc{BPV+Ilu(1e-4)} & 5 & 15.96 & & 233 & 292.91 \\
\textsc{BPV+Cholinc(1e-4)} & 5 & 2.00 & &  233 & 216.41\\
\textsc{BPV+Jacobi} & 5 & 10.71 & & 5 & 462.76 \\
\textsc{Gmres} & 4 & 6.62 & & 5 & 202.92 \\
\hline
\end{tabular}
\caption {Comparison of the computation time for N-S with $\nu=1$}
\end{center}
\end{table}

\newpage
\indent  From these tables, we can see that: (1) According to CPU
time, the preconditioner \textsc{Cholinc} is almost the same time as
the exact preconditioner, and both are better than \textsc{Jacobi}
and \textsc{Ilu} preconditioners. (2) The convergence rate of our
algorithm  is faster than \textsc{Bramble-Pasciak-Vassilev (BPV)}
algorithm in \cite{BramblePasciakVassilev_MC99}, and is comparable
to \textsc{Gmres} algorithm. (3) In our algorithm, $\tau_i$ can be
updated in each iteration, requiring no prior estimate on the
spectrum of Schur complement, but the convergence of  BPV algorithm
depends on the spectrum of preconditioners.

\indent The streamline of velocity and the contour of  pressure
using the algorithm in Section 3 are shown in figure 2, where we use
the adaptive Uzawa algorithm with the preconditioner obtained by the
incomplete Cholesky factorization with the tolerance of $10^{-1}$.
 \begin{figure}[!htbp]
\begin{minipage}[t]{0.5\linewidth}
  \centering
 \includegraphics [width=\textwidth]{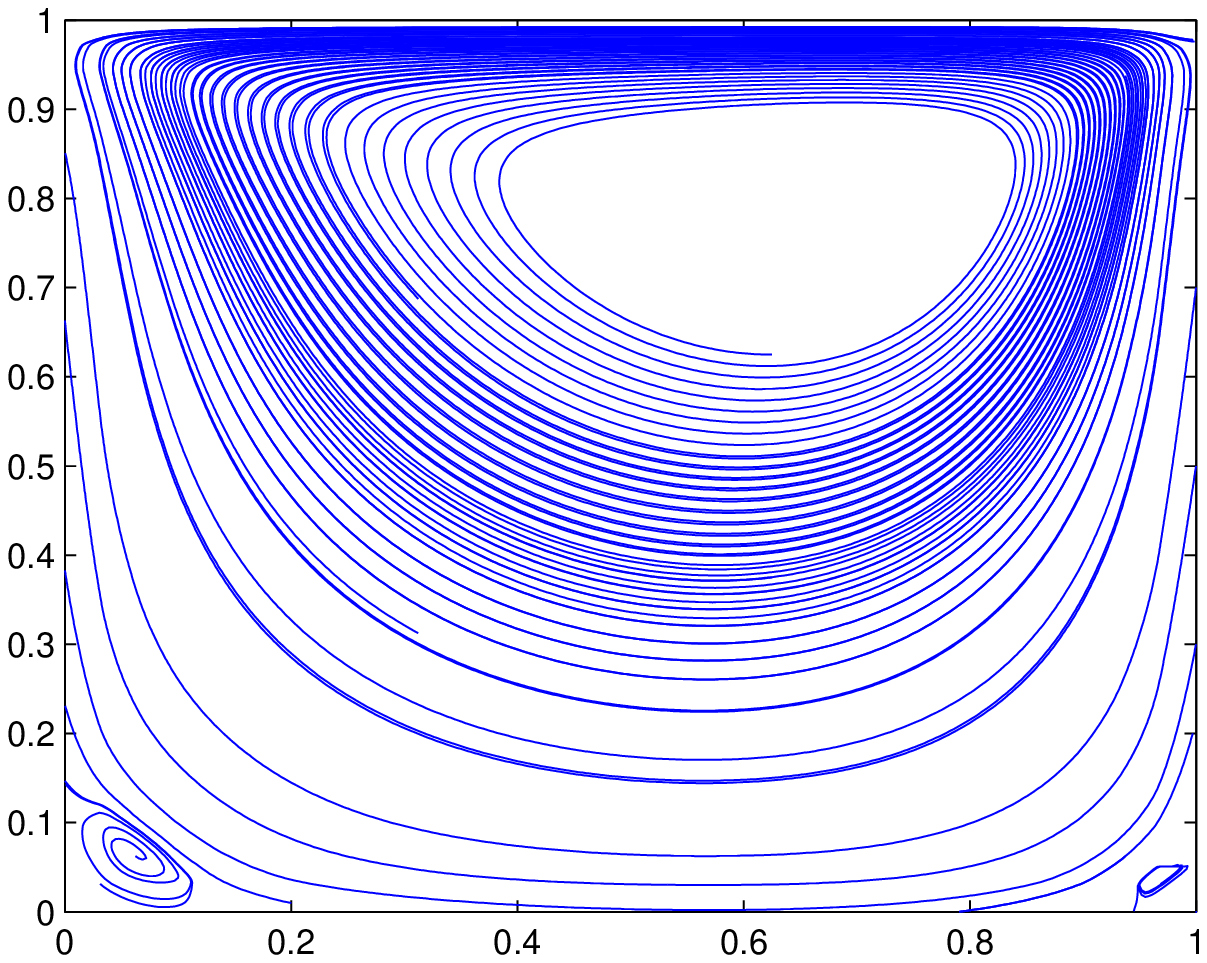}
\end{minipage}
 \begin{minipage}[t]{0.5\linewidth}
  \centering
 \includegraphics [width=\textwidth]{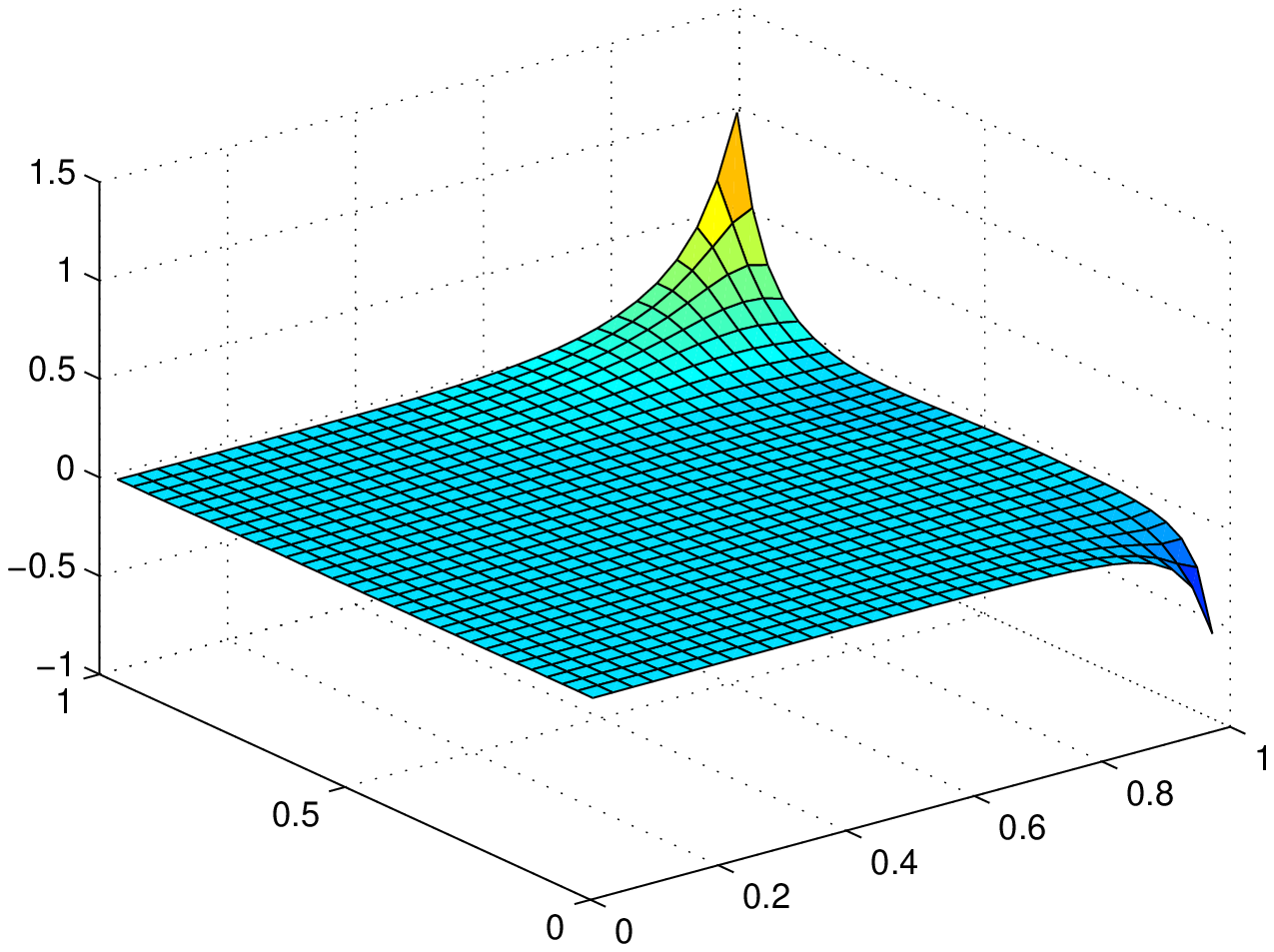}
\end{minipage}
\caption{Streamline (left) and pressure contour (right) of
Navier-Stokes problem ($32\times32$ grids, $\nu=0.01$)} \label{fig4}
\end{figure}
\section{Conclusion}

\indent In this paper, we investigate an adaptive Uzawa algorithm
with one variable relaxation parameter on generalized nonsymmetric
saddle point problems. Our work is closely related to the work in
\cite{BramblePasciakVassilev_MC99} and \cite{HuZou_MAA01}. In
\cite{HuZou_MAA01}, Hu and Zou introduced the adaptive Uzawa
algorithm with variable relaxation parameters for symmetric saddle
point problems with $D=0$, while we extend this algorithm to the
generalized nonsymmetric case. In
\cite{BramblePasciakVassilev_MC99}, Bramble \emph{et al.} discussed
Uzawa algorithm on nonsymmetric saddle point problems with $D=0$ and
proved its convergence under the assumption \eqref{eq1.11} and
\eqref{eq1.12}. We adopt their algorithm by adding a variable
relaxation parameter and prove its convergence without the condition
\eqref{eq1.11} and \eqref{eq1.12}, that is to say, without any prior
estimate on the spectrum of two preconditioned subsystems involved.
Our numerical experiments on Oseen problem and N-S problem
demonstrated the efficiency of our adaptive Uzawa algorithm. In
fact, in our computation experiments, when we choose  two variable
relaxation parameters just like in \cite{HuZou_MAA01}, we also can
get the convergence result. But its convergence theory is not
completed yet. Moreover,  the convergence result of nonlinear Uzawa
type algorithm with variable relaxation parameters also still needs
future
    work.

%\ack We would like to thank the editor and the referee for the
%valuable comments. % and helpful suggestions.
% H. Xiang is supported by the National Natural Science Foundation of China under grant 10901125 and  91130022.

%------------------------------------------------------------

\end{document}